\documentclass[12pt, final]{article}

\pdfoutput=1

\RequirePackage{amsfonts, amsmath} 

\RequirePackage{graphicx}                  
\RequirePackage{algorithm, algorithmic}    
                                           
\RequirePackage{color, colortbl}           
\RequirePackage{enumitem}

\usepackage{amsopn}
\RequirePackage{amsthm}
\RequirePackage{amsbsy}
\theoremstyle{plain} 
\newtheorem{theorem}{Theorem}[section]

\newtheorem{lemma}      [theorem]{Lemma}

\usepackage{multirow}
\usepackage{amsopn}
\usepackage{hyperref}
\usepackage{wrapfig}
\usepackage{boxedminipage}
\usepackage[normalem]{ulem}

\usepackage{pgfgantt}
\usepackage[font=small]{caption}

\definecolor{light-gray}{gray}{0.95}

\usepackage[backgroundcolor = light-gray, linecolor = gray,textsize = footnotesize]{todonotes}

\ifpdf
  \DeclareGraphicsExtensions{.eps,.pdf,.png,.jpg}
\else
  \DeclareGraphicsExtensions{.eps}
\fi

\newcommand{\Ac}{{\cal A}}

\newcommand{\Fc}{{\cal F}}

\newcommand{\arras}{\,\,\stackrel{{\rm a.s.}}{\longrightarrow}\,\,}

\newcommand{\Ex}{\mathbb{E}}

\newcommand{\I}{\mathbb{I}}

\newcommand{\Prob}{{\mathbb{P}}}
\newcommand{\Var}{{\rm Var}}

\def\zal2{z_{\alpha/2}}

\newcommand{\R}{{\mathbb{R}}}
\newcommand{\zero}{\boldsymbol{0}}

\RequirePackage{amsfonts, amsmath} 
\usepackage{fullpage}
\usepackage{tiaStyle}

\newcommand{\TheTitle}{Stochastic Newton and Quasi-Newton Methods for Large Linear Least-squares Problems}
\newcommand{\TheAuthors}{Chung, Chung, Slagel, and Tenorio}

\title{{\TheTitle}}

\author{
  Julianne Chung\thanks{Department of Mathematics and Computation Modeling and Data Analytics Division, Academy of Integrated Science, Virginia Tech, Blacksburg, VA
    (\href{mailto:jmchung@vt.edu}{jmchung@vt.edu}, \url{http://www.math.vt.edu/people/jmchung/}).}
  \and
  Matthias Chung\thanks{Department of Mathematics and Computation Modeling and Data Analytics Division, Academy of Integrated Science, Virginia Tech, Blacksburg, VA
    (\href{mailto:mcchung@vt.edu}{mcchung@vt.edu}, \url{http://www.math.vt.edu/people/mcchung/}).}
  \and
  J. Tanner Slagel\thanks{Department of Mathematics, Virginia Tech, Blacksburg, VA
    (\href{mailto:slagelj@vt.edu}{slagelj@vt.edu}, \url{https://www.math.vt.edu/people/slagelj/}).}
  \and
  Luis Tenorio\thanks{Department of Applied Mathematics and Statistics, Colorado School of Mines, Golden, CO
    (\href{mailto:ltenorio@mines.edu}{ltenorio@mines.edu}, \href{http://inside.mines.edu/\%7Eltenorio/}{http://inside.mines.edu/ltenorio/}).}
}

\ifpdf
\hypersetup{
  pdftitle={\TheTitle},
  pdfauthor={\TheAuthors}
}
\fi

\begin{document}

\maketitle
\begin{abstract}
  We describe stochastic Newton and stochastic quasi-Newton approaches to efficiently solve large linear least-squares problems
  where the very large data sets present a significant computational burden (e.g., the size may exceed computer memory or data are collected in real-time).
  In our proposed framework, stochasticity is introduced in two different frameworks as a means to overcome these computational limitations, and probability distributions that can exploit structure and/or sparsity are considered.
  Theoretical results on consistency of the approximations for both the stochastic Newton and the stochastic quasi-Newton methods are provided. The results show, in particular, that stochastic Newton iterates, in contrast to stochastic quasi-Newton iterates, may not converge to the desired least-squares solution. Numerical examples, including an example from extreme learning machines, demonstrate the potential applications of these methods.

\end{abstract}

{\bf Keywords:} stochastic approximation, stochastic Newton, stochastic quasi-Newton, least-squares, extreme learning machine.




\section{Introduction} \label{sec:introduction} 

In this paper, we are interested in linear problems of the form
\begin{equation}
	\label{eq:inverseprob}
	\bfb = \bfA \bfx_\true + \bfepsilon,
\end{equation}
where $\bfb \in \bbR^m$ contains the observed data, $\bfx_\true \in \bbR^n$ is the desired solution, $\bfA \in \bbR^{m \times n}, m \gg n$, has full column rank, and $\bfepsilon$ representes noise modeled as random. Here, we assume that $\bfepsilon$ has zero mean and covariance matrix $\sigma^2\bfI_m$, where $\bfI_m$ is the $m \times m$ identity matrix.
The goal is to compute the unique least-squares (LS) solution,
\begin{align}
\label{eq:argmin}
\widehat\bfx = \argmin_{\bfx \in \mathbb{R}^n} \ f(\bfx) = \tfrac{1}{2}\norm[]{\bfA\bfx-\bfb}^2,
\end{align}
where $\norm[]{\,\cdot\,}$ denotes the two norm. Typical LS solvers such as iterative methods for large, sparse problems may require many multiplications with $\bfA$ and its transpose \cite{Nocedal2006,hansen2012least}.  This can be prohibitively expensive for many problems of interest.
We overcome these limitations by first reformulating the LS problem as a stochastic optimization problem as follows.
Let $\bfW\in \bbR^{m\times \ell}, m > \ell$, be a random matrix with
\[
\bbE \left(\, \bfW \bfW\t\,\right)  = \beta \bfI_m, 
\]
where $\beta > 0$ and $\bbE$ denotes the expectation operator. Then
\begin{gather*}
\bbE \,\norm[]{\bfW\t\left(\bfA\bfx - \bfb\right)}^2 \,
 = (\bfA\bfx - \bfb)\t \bbE \left(\,\bfW \bfW\t\,\right) (\bfA\bfx - \bfb) = \beta \norm[]{\bfA\bfx - \bfb}^2.
\end{gather*}
Since the solution in~\eqref{eq:argmin} is \emph{equivalent} to the solution of the following stochastic optimization problem,
\begin{align}
	\label{eq:stochmin}
	\min_{\bfx} \ \bbE f_\bfW(\bfx),\qquad \mbox{where } \quad f_\bfW(\bfx) = \tfrac{1}{2\beta} \norm[]{\bfW\t\left(\bfA\bfx - \bfb\right)}^2,
\end{align}
the goal of this work is to solve~\eqref{eq:stochmin}.
 Two approaches have been used in the literature to approximate~\eqref{eq:stochmin}: \emph{Sample Average Approximation} (SAA) and \emph{Stochastic Approximation} (SA) \cite{robbins1951,shapiro2014lectures}.  In SAA methods, the goal is to solve a sampled LS problem, where one major challenge is to find a good sampling matrix $\bfW$ (e.g., using randomization techniques).
On the other hand, SA is an iterative approach where different realizations of the sampling matrix $\bfW$ are used at each iteration to update the approximate solution.
In this work we will only consider the SA approach.
%

Our main contributions are summarized as follows. We reformulate the problem as a stochastic optimization problem, as described above, and develop stochastic Newton and stochastic quasi-Newton methods for solving large LS problems. By introducing stochasticity and allowing sparse realizations of $\bfW$, we can handle large problems where only pieces of data are readily available or in memory at a given time. We prove almost sure convergence of stochastic quasi-Newton estimators to the desired LS solution.
%
Furthermore, we show that for systems that are inconsistent but whose coefficient matrix has full column rank, stochastic Newton iterates may not converge to the unique LS solution~$\widehat{\bfx}$, but instead will converge to the solution of a different (not necessarily nearby) problem (see Section~\ref{sec:theoretical_results}). By exploiting updates for the inverse Hessian approximation, our stochastic quasi-Newton method can produce good approximations and be computationally more efficient than standard LS methods, especially for very large problems.

A similar stochastic reformulation was described in \cite{Le2016,drineas2011faster,woodruff2014sketching}, but a fundamental difference with our work is that these approaches use an SAA, rather than an SA method to approximate $\widehat \bfx$.  SA methods have been intensively studied, especially in the last few years.  Although much of the literature has focused on general objective functions, general convergence theories, and gradient-based methods \cite{Kushner1997, Bottou1998,benveniste2012,needell2016stochastic,strohmer2009randomized}, only recently has the focus shifted to higher-order methods (e.g., stochastic Newton) and efficient update methods that can incorporate curvature information (e.g., stochastic quasi-Newton) \cite{Bottou2016,Xu2016}.  Repeated sampling in a higher order method was considered in \cite{pilanci2016iterative,Haber2012}.  However, \cite{Haber2012} does not consider sparse sampling, and the algorithm described in \cite{pilanci2016iterative} requires the full gradient of the objective function, which we do not require here, and uses only curvature information from the current sample, whereas our stochastic quasi-Newton approach can integrate all previous samples in an efficient update scheme.  Gower and Richt{\'a}rik proved in \cite{Gower2015} that for consistent linear systems (e.g.,~\eqref{eq:inverseprob} with $\bfepsilon = \bfzero$), the stochastic Newton method with unit step size converges to the desired solution.  However, they do not consider LS problems.  Furthermore, our work is different from recent works on stochastic quasi-Newton methods that focus on computing matrix inverses~\cite{gower2016} and solving empirical risk minimization problems~\cite{gower2016s}.

An outline for this paper is as follows.  In Section~\ref{sec:overview}, we describe the  stochastic Newton and stochastic quasi-Newton method for solving LS problems whose
consistency results are presented and discussed in Section~\ref{sec:theoretical_results}.
Numerical results in Section~\ref{sec:numerics} demonstrate the potential of our approaches, and conclusions are provided in Section~\ref{sec:conclusions}. Consistency proofs are provided in Section~\ref{sec:appendix}.

%

%


\section{Stochastic Newton and quasi-Newton} 
\label{sec:overview}
In this section, we describe SA methods for computing a solution to~\eqref{eq:stochmin}.  Given an initial vector $\bfx_{0} \in \mathbb{R}^{n}$, SA methods define a sequence of iterates
\begin{align} \label{eq:alg}
	\bfx_{k}= \bfx_{k-1}+ \alpha_{k}\bfs_{k},
\end{align}
where $\{\alpha_{k}\}$ is a sequence of step sizes and $\{\bfs_{k}\}$ is a sequence of search directions such that $\bfs_k$ depends on the iterate $\bfx_{k-1}$ and the random variables $\bfW_{1},\dots,\bfW_{k}$.
The most common SA approach is the \emph{stochastic gradient method}, where $\bfs_k = -\nabla f_{\bfW_k}(\bfx_{k-1})$ is the sample gradient,
\begin{equation} \label{eq:samplegradient}
	\nabla f_\bfW(\bfx) = (\bfW\t\bfA)\t(\bfW\t\bfA\bfx-\bfW\t\bfb)=  \bfA^\top\bfW\bfW^\top(\,\bfA\bfx - \bfb\,)\,,
\end{equation}
evaluated at $\bfx_{k-1}$.
The popularity of the stochastic gradient method stems from its proven consistency properties and its easy implementation. However, the stochastic gradient method is known to converge slowly \cite{Xu2016},  thus higher order methods are desired and discussed below.

\medskip
For the \emph{stochastic Newton} method, the search direction is typically defined as
	\begin{equation} \label{eq:Newtonstep}
		\bfs_{k}  = - \left(\nabla^{2} f_{\bfW_{k}}\right)^{\dagger}\nabla f_{\bfW_{k}}(\bfx_{k-1}),
	\end{equation}
	where the sample Hessian is given by
		$\nabla^{2} f_\bfW = \bfA\t \bfW \bfW\t\bfA \,,$
	 and $\dagger$ denotes the Moore-Penrose pseudoinverse.
	Using properties of the pseudoinverse,~\eqref{eq:Newtonstep} can be reduced to
	\begin{align} \bfs_{k} = -(\bfW_{k}\t\bfA)^{\dagger}\bfW_{k}\t(\,\bfA\bfx_{k-1}-\bfb\,).
	\label{eq:pseudo3}
	\end{align}
	For the stochastic Newton method, we use a step size $\alpha_k$ that
	satisfies the following frequently used conditions
 \begin{equation} \label{eq:steplength}
 \sum_{k=1}^{\infty}\alpha_{k} = \infty \qquad \text{and}  \qquad \sum_{k=1}^{\infty}\alpha_{k}^{2} < \infty\,.
 \end{equation}
	%
For the \emph{stochastic quasi-Newton method}, we use the search direction given by
	\begin{align}
	\bfs_{k}  & = \label{eq:StochQuas} -\bfB_{k}\nabla f_{\bfW_{k}} (\bfx_{k-1}),
	\end{align}
where $\{\bfB_{k}\}$ is a sequence of \emph{random} symmetric positive definite (SPD) matrices. For the convergence results proved in Section~\ref{sec:theoretical_results}, we require
that $\lambda_{\text{max}}(\bfB_{k}) \leq B$ almost surely (a.s.) for some constant $B >0$ and for all $k$, with
	\begin{equation} \label{eq:scalarc}
		\sum_{k=1}^{\infty}\alpha_{k}\,\lambda_{\text{min}}\left(\bfB_{k}\right) = \infty\, \,\,\mathrm{a.s.} \qquad \text{and} \qquad \sum_{k=1}^{\infty}\alpha_{k}^{2} < \infty \,,
	\end{equation}
	where $\lambda_{\text{min}}\left(\bfB_{k}\right)$ and $\lambda_{\text{max}}\left(\bfB_{k}\right)$ denote the minimal and maximal eigenvalues of $\bfB_{k}$, respectively.
	Note that the special case where $\bfB_{k}=\bfI_n$ is nothing more than the stochastic gradient method.
	The sequence $\{\bfB_{k}\}$ is assumed to be independent of $\bfW_{k}$, but may depend on  $\bfW_{1}, \dots, \bfW_{k-1}$.

These two algorithms provide different ways to estimate $\widehat{\bfx}$. The choice of $\{\alpha_k\}$ and the distribution of $\bfW$ fully determine the stochastic Newton method, and these along with the choice of $\{\bfB_k\}$ fully determine the stochastic quasi-Newton method.  Below we briefly summarize choices for $\{\alpha_k\}$, the distribution of $\bfW$, and the sequence $\{\bfB_k\}$, with particular attention on  the conditions required to prove consistency of the approximations.



\smallskip
\paragraph{\rm \textbf{Choice of $\{\alpha_k\}$}}
Selecting a good step size $\alpha_k$ (or learning rate, as it is referred to in machine learning) is critical. 
A variety of methods have been proposed to improve convergence rates, see for instance \cite{Bottou1998,Spall2005,Diniz1997}. In this paper we restrict ourselves to step size choices that comply with conditions~\eqref{eq:steplength} and~\eqref{eq:scalarc} for stochastic Newton and stochastic quasi-Newton, respectively. In our numerical experiments in Section~\ref{sec:numerics}, we use the harmonic sequence $\alpha_k = {1}/{k}$, see~\cite{Robbins1985}.


\medskip
\paragraph{\rm\textbf{Distributional choice for} $\bfW$}
As discussed in the introduction, the solutions to problems~\eqref{eq:argmin} and~\eqref{eq:stochmin} are equivalent if $\bbE( \bfW\bfW\t ) = \beta \bfI_m$. Extensions to include a general (known) positive definite covariance matrix $\bfGamma$ would require adjusting the sampling matrix such that $\bbE(\bfW\bfW\t)= \bfGamma^{-1}$. Including a positive diagonal covariance matrix would simply require scaling the rows of $\bfA$ (i.e., solving a weighted LS problem).  However, for simplicity of presentation, we consider $\bfGamma = \sigma^2\, \bfI_m$.

%
\medskip
\noindent Three distributional choices for $\bfW$ are considered:
\begin{enumerate}[leftmargin = 3ex,itemsep = 1ex]
\item \emph{Random sparse matrices.}
Let $\bfW \in \bbR^{m \times \ell}$ be a random matrix with i.i.d. random elements $w_{ij}$ where, for a fixed $0<\psi\leq 1$,
$w_{ij}$ takes the values $\pm \sqrt{\beta/\ell\psi}$ each with probability $\psi/2$ and the value zero with probability $1-\psi$.
				%
It is straightforward to verify that $\bbE(\, \bfW\bfW\t\,) = \beta \bfI_m$. Notice that as $\psi$ gets closer to zero, more sparsity is introduced in $\bfW$. It is worth mentioning  that this choice of $\bfW$ is a generalization of Achlioptas random matrix ($\psi = 1/3$ and $\beta = \ell$) and the Rademacher distribution ($\psi = 1$ and $\beta = \ell$), see~\cite{Achlioptas2003,Haber2012}.\label{pg:sparserandom}
\item \emph{Generalized Kaczmarz matrices.}
For $i = 1,\ldots,p$, let $\bfQ_i \in \bbR^{m \times \ell_i}$ be such that  $\bfQ = [\bfQ_1,\ldots,\bfQ_p] \in \bbR^{m \times m}$ is an orthogonal matrix. Define the distribution of $\bfW$ to be
uniform on $\{\bfQ_1,\ldots,\bfQ_p\}$. Then
\begin{equation*}
\bbE \left(\bfW\bfW\t \right) = \tfrac{1}{p} \sum_{i = 1}^p \bfQ_i \bfQ_i\t = \tfrac{1}{p} \,\bfI_m.
\end{equation*}
 Notice that selecting $\bfQ = \bfI_m$, or any permutation of it, together with the stochastic Newton direction~\eqref{eq:Newtonstep} leads to the well-known randomized Kaczmarz method if $\ell_i = 1$ for all $i$, and to the randomized block Kaczmarz method otherwise \cite{kaczmarz1937angeniherte, Needell2014,Gower2015}.
Choosing the elements of $\bfQ$ to be $\pm 1$ (or in $\{0,\pm 1\}$) leads to (sparse) randomized Hadamard matrices \cite{Hedayat1978,Boutsidis2012}. Notice, that sparsity may be introduced by the particular choice of $\bfQ$ and that the number of columns in the $\bfQ_i$'s can differ.\label{pg:kaczmarzW}
			\item \emph{Sparse Rademacher matrices.}
				Fix $p \leq m$. The columns $\bfw_i$ of $\bfW \in \mathbb{R}^{m \times \ell}$  are i.i.d. and each column can be any $m\times 1$ vector with $p$-nonzero entries in $\{\pm 1\}$ with equal probability. Hence, conditional on a vector configuration, $C$, of $p$ ones and $m-p$ zeros, each column $\bfw_i$ has
conditional expectation
				\[
					\Ex(\,\bfw_i\bfw_i\t\mid C\,) = \bfI_{m,C},
				\]
				where $ \bfI_{m,C}$ is the diagonal matrix with the configuration $C$ in the diagonal. It follows that
				\[
					\Ex(\,\bfw_i\bfw_i\t\,) = \Ex\,\Ex(\,\bfw_i\bfw_i\t\mid C\,)=\frac{
					\begin{pmatrix} m-1\\ p-1\end{pmatrix}}{
					\begin{pmatrix} m\\ p\end{pmatrix}
					}\,\bfI_m = \tfrac{p}{m}\,\bfI_m,
				\]
				and therefore $\Ex(\,\bfW\bfW\t\,) = (\ell\,p/m)\,\bfI_m$.
%
				Note  that the case $p = m$ generates full Rademacher matrices. The distinction between the other choices of $\bfW$ is that entries of the sparse Rademacher matrices are not i.i.d., as in the random sparse matrices, and do not necessarily come from partitions of orthogonal matrices, as in the generalized Kaczmarz matrices.
		\end{enumerate}

\medskip
\paragraph{\rm \textbf{Choice of} $\{\bfB_{k}\}$}
The matrices $\bfB_{k}$ should approximate the inverse Hessian $(\bfA\t\bfA)^{-1}$ and should be SPD.
Furthermore, the
convergence results require some control of the smallest and largest eigenvalues of $\bfB_{k}$. Thus, we define a sequence of SPD matrices $\bfB_k$ that converges a.s.\ to
a small perturbation of $(\bfA\t\bfA)^{-1}$ and has appropriate behavior of the smallest and largest eigenvalues. Since for any $\lambda_1>0$
\[
\tfrac{\lambda_1}{k}\bfI_n + \tfrac{1}{k}\sum_{i=1}^{k}\bfA\t\bfW_{i}\bfW_{i}\t\bfA\,\, \arras\,\, \bfA\t\bfA,
\]
it follows that for any $\lambda_2\geq 0$,
\begin{eqnarray*}
\widetilde{\bfB}_k = \lambda_{2}\,\bfI_{n} + k \left(\,\lambda_{1}\bfI_{n}+ \sum_{i=1}^{k}\bfA\t\bfW_{i}\bfW_{i}\t\bfA\,\right)^{-1} &\arras\,\, & \lambda_{2}\,\bfI_{n} + (\bfA\t\bfA)^{-1},
\end{eqnarray*}
and $\lambda_\mathrm{max}(\widetilde{\bfB}_k) \arras  \lambda_2 + \lambda_\mathrm{max}(\, (\bfA\t\bfA)^{-1}\,)$. Hence, for a fixed $B> \lambda_2 + \lambda_\mathrm{max}(\, (\bfA\t\bfA)^{-1}\,)$,
we can define
\begin{equation} \label{eq:Bupdate}
\bfB_{k} =
			\begin{cases}
				\widetilde{\bfB}_k,  & \text{if } \lambda_{\text{max}}(\widetilde{\bfB}_k) \leq B,  \\
				\bfB_{k-1}, & \text{else,}\\
			\end{cases}
\end{equation}
with $\bfB_0=\widetilde{\bfB}_0$ an arbitrary SPD matrix such that $\lambda_{\text{max}}\left(\bfB_{0}\right) \leq B$. The Woodbury formula can be used to efficiently update $\bfB_k$, see Section~\ref{sec:numerics} for details. Notice that  $\widetilde \bfB_k$ is a sample approximation of $(\bfA\t\bfA)^{-1}$, which is similar to sample approximations that have been studied for covariance/precision matrix approximations, and that the regularization term $\lambda_1\, \bfI_n$ can be replaced by other appropriate SPD matrices \cite{fan2016overview}.




%

\section{Theoretical results}
\label{sec:theoretical_results}
In this section we study consistency of the stochastic quasi-Newton and stochastic Newton method.
We start by briefly reviewing some notation and definitions from probability theory.

\smallskip
Let $(\Omega,\Ac,\Prob)$ be a probability space, and let $\{\calF_{k}\}$ be a sequence of sub-$\sigma$-algebras of $\Ac$. Then
$\{\calF_{k}\}$ is a called a filtration if $\calF_{k} \subset \calF_{k+1}$ for all $k \in \mathbb{N}$. A sequence of random variables $\{u_{k}\}$ on $(\Omega,\Ac)$ is
said to be \textit{adapted to} $\{\calF_{k}\}$ if $u_{k}$ is $\calF_{k}$-measurable for all $k \in \mathbb{N}$. The $\sigma$-algebra generated by
the random variables $\{u_{i}\}_{i=1}^{k}$ is denoted by $\sigma\left(u_{i} : i < k\right)$. We use $\I_A$ to denote the indicator
function of the set $A$. Our convergence results make use of the following quasimartingale convergence theorem. For a proof and further details on quasimartingales, we refer the interested reader to \cite{Fisk1963,Metivier1982}.

\begin{theorem}[Quasimartingale convergence theorem] \label{thm:quasi}
Let $\{X_n\}$ be a sequence of non-negative random variables adapted to a filtration $\{\calF_{k}\}$ and such that $\Ex\,X_n<\infty$ for all $n$.
Let $Z_k = \bbE \left(X_{k+1}-X_{k} | \calF_{k} \right)$. Then, if
  	\begin{equation*}
 		\sum_{k=0}^{\infty} \bbE \left[\I_{Z_k\geq 0}\,(X_{k+1}-X_{k})\right] < \infty,
 	\end{equation*}
then there is an integrable, non-negative random variable $X$ such that $X_n \arras X$.
\end{theorem}

Now we are ready to state the convergence of the stochastic quasi-Newton to the LS solution. A proof can be found in the appendix (Section~\ref{sec:appendixSQN}).
\begin{theorem}[Stochastic quasi-Newton method] \label{thm:SQN}
	Let  $\bfA \in\R^{m\times n}$ have rank $n$ and $\bfb \in\R^m$. Let $\{\bfW_k\}$ be a sequence of i.i.d. $m\times\ell$ random matrices with $\Ex\left(\bfW_k\bfW_k^\top\right) = \beta\bfI_{m}$ and $\Ex(\,\norm[]{\bfW_k\bfW_k^\top}^2\,)<\infty$. Define $\Fc_k = \sigma(\,\bfW_i;\,i<k\,)$. Let $\{\alpha_k\}$ be a positive sequence of scalars, and $\{\bfB_k\}$ be a sequence of random $n\times n$ SPD matrices adapted to $\left\{\Fc_k \right\}$ such that for some $B>0$
\begin{equation*}
 \|\bfB_k\|=\lambda_\mathrm{max}(\bfB_k)\leq B\qquad \mathrm{a.s.} \quad \forall k,
\end{equation*}
and
	\begin{equation*}
		\sum_{k=1}^{\infty} \alpha_k\,\lambda_{\min}(\bfB_k) =\infty\,\,\mathrm{a.s.}\quad\mbox{and}\quad \sum_{k=1}^{\infty} \alpha_k^2  < \infty.
	\end{equation*}
	Set $\widehat{\bfx}=(\bfA^\top\bfA)^{-1}\bfA^\top\bfb$, $\widehat{\bfb} = \bfA\widehat{\bfx}$, and let $\bfx_0\in \R^n$ be an arbitrary initial vector. Define
         \begin{eqnarray*}
	\bfx_{k} &= & \bfx_{k-1} + \alpha_k\bfs_k,\\
	\bfs_k     &= &-\bfB_k\bfA^\top\bfW_k\bfW_k^\top(\bfA\bfx_{k-1} - \bfb), \mbox{ and}\\
        \bfb_k & = & \bfA\bfx_{k-1}.
\end{eqnarray*}
 Then: 
	\begin{enumerate}[label=(\roman*),leftmargin = 5ex]
		\item The matrix $\bfC=\Ex(\,\bfW_k\bfW_k^\top\bfW_k\bfW_k^\top\,)$ is defined and SPD.
		\item If $e_k = \|\,\bfb_k - \widehat{\bfb}\,\|^2$, then $\Ex \,e_k <\infty$ and $\Ex (\,\|\,\bfs_k\,\|^2\,)<\infty$ for all $k$.
		\item $\bfx_k\arras \widehat{\bfx}$.
	\end{enumerate}
\end{theorem}

\medskip

The following result shows that the stochastic Newton method does not necessarily converge to the LS solution $\widehat \bfx$. For a proof see Section~\ref{sec:appendixSN}
in the appendix.

\begin{theorem}[Stochastic Newton method]\label{thm:SN}
	Let  $\bfA \in\R^{m\times n}$ have rank $n$ and $\bfb \in\R^m$. Let $\{\bfW_k\}$ be a sequence of i.i.d. $m\times\ell$ random matrices with $\Ex(\,\bfW_k\bfW_k^\top\,) = \beta\bfI_{m}$. Let $\bfH_k =(\bfW_k\t\bfA)^{\dagger}\bfW_k\t$ and assume that $\bfC=\bbE(\,\bfH_k\t\bfH_k \,)$ is defined and finite. Let $\{\alpha_k\}$ be a positive sequence of scalars such that
	\begin{equation*}
  	\sum_{k=1}^{\infty} \alpha_k  =\infty \quad\mbox{and}\quad \sum_{k=1}^{\infty} \alpha_k^2  < \infty.
	\end{equation*}
	Set $\bfP = \Ex\,\bfH_k$, $\widetilde{\bfx}=(\bfP\bfA)^{-1}\bfP \bfb$, and let $\bfx_0\in \R^n$ be an arbitrary initial vector. Define
	\begin{eqnarray*}
		\bfx_{k} &= & \bfx_{k-1} + \alpha_k\bfs_k\\
		\bfs_k &= &-\bfH_k\left(\,\bfA\bfx_{k-1}-\bfb\,\right).
	\end{eqnarray*}
 	Then: 
	\begin{enumerate}[label=(\roman*),leftmargin = 5ex]
		\item The matrices $\bfC$ and $\bfA\bfC\bfA\t$ are symmetric positive semi-definite, and $\bfP\bfA$ is SPD.
		\item If $e_k = \|\,\bfx_k - \widetilde{\bfx}\,\|^2$, then $\Ex \,e_k <\infty$ and $\Ex (\,\|\,\bfs_k\,\|^2\,)<\infty$ for all $k$\,.
		\item $\bfx_k\arras \widetilde{\bfx}$\,.
	\end{enumerate}
\end{theorem}

Although the number of iterations can not be too large in practice,
the consistency results above for stochastic Newton and stochastic quasi-Newton are important for analyzing these new methods. In the rest of this section, we consider the feasibility of the convergence criteria for both methods and then provide some insight into the discrepancy between the stochastic Newton estimator $\widetilde \bfx$ and the desired LS solution.

\medskip
In addition to $\mathbb{E}(\bfW \bfW\t ) = \beta\bfI_m,$ Theorem~\ref{thm:SQN} requires $\mathbb{E} (\,\norm[]{\bfW_{k}\bfW_{k}\t}^2\,)$ to be finite, while Theorem~\ref{thm:SN} requires $\bfP$
to be finite.  However, these additional assumptions are not too restrictive.  In particular, for all of the choices of $\bfW$ described in Section \ref{sec:overview}, these expectations are trivially finite because the entries of $\bfW$ take only finitely many values.
As for the choice of $\bfB_k$ in the stochastic quasi-Newton method, the assumption that $\bfB_k$ and $\alpha_k$ satisfy~\eqref{eq:scalarc} is not difficult to meet. For example,
when $\lambda_{2}>0$
the choice of $\bfB_{k}$ in  Section \ref{sec:overview} satisfies these conditions trivially because
\[
\lambda_{2} \leq \lambda_{\text{min}} (\bfB_{k} ) \leq  \lambda_{\text{max}} (\bfB_k )\leq B.
\]
On the other hand, $\lambda_\mathrm{max}(\bfW)$ is bounded for the choices of $\bfW$ in this paper, and therefore there is a $C>0$ such that
$\|\bfA\t\bfW_i\bfW_i\t\bfA\| \leq C$ for all $i$ and
\[
\left\|\lambda_1\,\bfI_n + \sum_{i=1}^k\bfA\t\bfW_i\bfW_i\top \bfA\right\| \leq \lambda_1 + k C,
\]
which implies that when $\lambda_2=0$ we have
\[
\lambda_\mathrm{min}(\bfB_k) \geq \frac{1}{\lambda_1/k + C} \geq \frac{1}{\lambda_1 + C}>0
\]
and therefore conditions~\eqref{eq:scalarc} are satisfied.

\begin{figure}[bthp]
  \hspace*{-5ex}\includegraphics[width=1.1\textwidth]{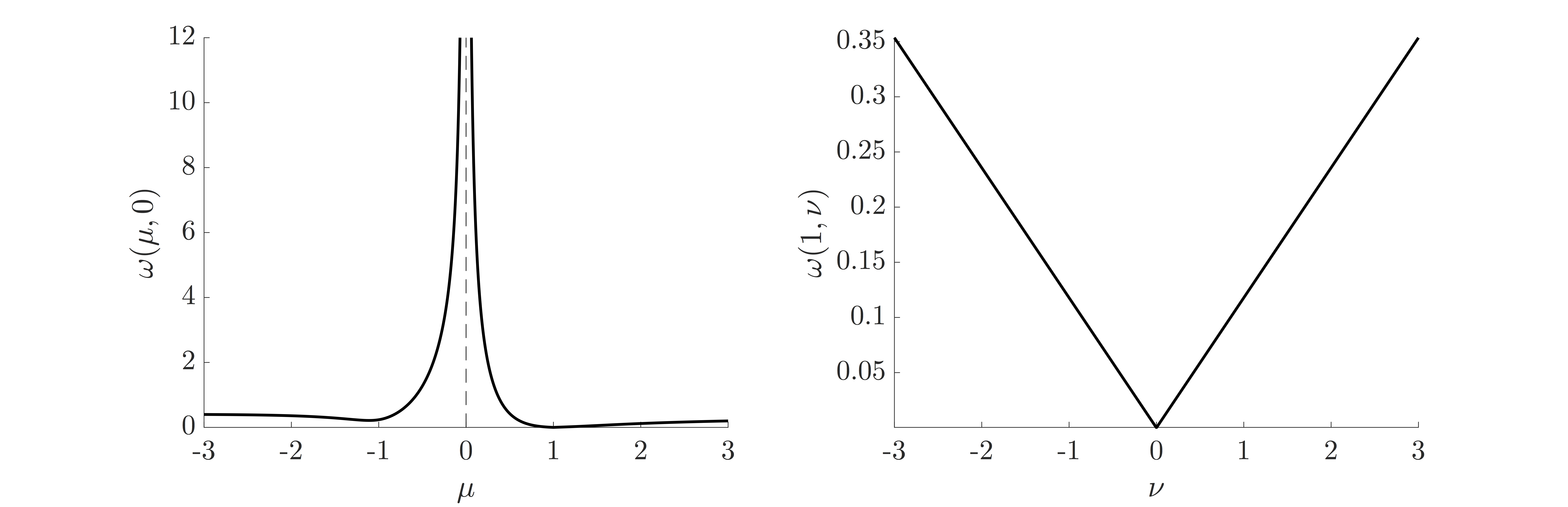}
  \caption{Error $\omega(\mu,\nu)$ of the stochastic Newton solution $\widetilde \bfx$ compared to $\widehat \bfx$. In the plot on the left, $\nu=10$ and we vary $\mu.$  Notice that a pole exists at $\mu = 0$, where the relative error becomes arbitrarily large. The plot on the right illustrates the impact of varying $\nu$ for fixed $\mu=1$.}\label{fig:Example0}
\end{figure}

\medskip
Next, we provide some insight regarding the potential discrepancy between the desired LS solution $\widehat\bfx$ and the solution to which the stochastic quasi-Newton method converges, namely,
	\begin{equation}\label{eq:tilde}
\widetilde{\bfx}  = (\bfP\bfA)^{-1}\bfP\bfb  = (\,\mathbb{E} \,[(\bfW\t\bfA)^{\dagger}\bfW\t\bfA\,]\,)^{-1}\,\mathbb{E}[\,(\bfW\t\bfA)^{\dagger}\bfW\t\,]\,\bfb .
\end{equation}
The difference between $\widehat\bfx$ and $\widetilde \bfx$ depends on $\bfP$, but we can say the following. Assuming that the noise has zero mean and covariance matrix
$\Var( \bfepsilon) =\sigma^2 \bfI_m$, we have
\begin{eqnarray*}
\Ex\,\widehat \bfx & = & \bfx_{\true},\qquad \Var(\widehat \bfx)  =  \sigma^2\,(\bfA\t\bfA)^{-1},\\
\Ex\,\widetilde{\bfx} & = & \bfx_{\true},\qquad \Var(\widetilde{\bfx})  = \sigma^2\, (\bfP\bfA)^{-1}\bfP\bfP\t(\bfA\t\bfP\t)^{-1}.
\end{eqnarray*}
This shows that $\widehat{\bfx}$ and $\widetilde{\bfx}$ are both unbiased estimators of $\bfx_{\true}$, but by the Gauss-Markov theorem, $\widehat{\bfx}$ is
expected to have smaller variance.
Consider the following simple example:

\smallskip
\paragraph{Example}
Consider the LS problem, where
\begin{equation*}
  \bfA =
  \begin{bmatrix} \mu & 0 \\ 0 & 1 \\ 1 & -1 \end{bmatrix}
  \qquad \text{and}  \qquad
  \bfb = \begin{bmatrix} 1 \\ 1 \\ \nu \end{bmatrix},
\end{equation*}
for some fixed $\mu,\nu\in\bbR$. We compare the LS solution and the solution obtained via stochastic Newton with random Kaczmarz vectors $\bfw\in\bbR^{m\times 1}$ (see page~\pageref{pg:kaczmarzW}). It is easy to see that in this case we have
\[
\bfP = \bfA\t\bfH\quad \mbox{with}\quad \bfH = \mathrm{diag}\{1/\|\bfa_1\|^2,1/\|\bfa_2\|^2,1/\|\bfa_3\|^2\},
\]
where $\bfa_i$ are the rows of $\bfA$.
It follows that $\widetilde{\bfx}$ minimizes the weighted LS functional
\(
(\bfb-\bfA\bfx)\t\bfH(\bfb-\bfA\bfx).
\)
We obtain the following solutions:
\begin{equation*}
  \widehat\bfx = \frac{1}{2\mu^2 + 1}
  \begin{bmatrix}
    2\mu + \nu + 1 \\[1ex] \mu - \mu^2\nu + \mu^2 + 1
  \end{bmatrix}
  \qquad \text{and}  \qquad
  \tilde\bfx = \frac{1}{4}
  \begin{bmatrix}
    1 + \nu + {3}/{\mu} \\[1ex]
    3 - \nu + {1}/{\mu}
  \end{bmatrix},
\end{equation*}
respectively. The covariance matrices of $\widehat{\bfx}$ and $\widetilde{\bfx}$ are
\[
\Var(\widehat{\bfx}) = \frac{\sigma^2}{2\mu^2 + 1}\begin{bmatrix}
2 & 1\\
1 & \mu^2 + 1
\end{bmatrix},\qquad
\Var(\widetilde{\bfx}) = \frac{\sigma^2}{4\mu^2}\begin{bmatrix}
2\mu^2 + 9 & 8\mu^2 + 3\\
8\mu^2 + 3 & 10\mu^2 + 1
\end{bmatrix}.
\]
It is clear that the variances of the components $\widetilde{\bfx}$ can be much larger than those of  $\widehat{\bfx}$. The solution $\widetilde{\bfx}$ would have
smaller variance if the covariance matrix of the noise was proportional to $\bfH^{-1}$ instead of $\bfI_m$.
 Figure~\ref{fig:Example0} shows the error $\omega(\mu,\nu) = \|\widehat\bfx - \widetilde\bfx\|$ for various choices of $\mu$ and $\nu$. The left panel
shows that  $\omega\to \infty$ as $\mu \to 0$, which makes sense as the first row of $\bfA$ becomes all zeros.
The right panel shows that even for $\mu \neq 0,$ a significant error can be incurred by varying $\nu$ -- and therefore the ``observation vector'' $\bfb$.

\medskip
Although the difference between $\widetilde \bfx$ and $\widehat \bfx$ can be significant, there are cases where they are identical.
Some previous works have studied the problem of how close $\widetilde \bfx$ is to $\widehat \bfx$, e.g., \cite{woodruff2014sketching,drineas2011faster}.  However, their assumptions do not apply to our matrix $\bfP.$
%
For our problem, $\widetilde \bfx = \widehat \bfx$ when the linear system is consistent, 
since in this case, $\bfP\bfA\widehat{\bfx}=\bfP\bfb$, or when
$ ( \bfA\t \bfA)^{-1}\bfA\t = ( \bfP\bfA )^{-1}\bfP$, which is equivalent to $\text{null}(\bfA\t) \subseteq \text{null}(\bfP)$, since $\bfA\widehat{\bfx}-\bfb \in \text{null}(\bfA\t)$ implies that $\bfP(\bfA\widehat{\bfx}-\bfb)=\bf0$.  In the example above, this occurs when $\nu={1}/{\mu}-1$ (e.g., $\mu=1$ and $\nu=0$).

\section{Numerical experiments} \label{sec:numerics}
In this section we showcase the numerical performance of the discussed methods. A general framework for our proposed methods is summarized in Algorithm~\ref{alg:SQN}, where the main distinction among the stochastic gradient, Newton, and quasi-Newton methods is the choice of $\bfB_k$ (see line~\ref{alg:update}). For the stochastic gradient method, $\bfB_k = \bfI_n$, and for the stochastic Newton method, $\bfB_k =  \left(\bfA\t\bfW_k\bfW_k\t\bfA \right)^{\dagger}$.  For the stochastic quasi-Newton method, we will use the choices of $\bfB_k$ described in Section~\ref{sec:overview} with $\lambda_{1}>0$ and $\lambda_{2}=0$.
To avoid large inversions when computing $\bfB_k$, we use the Woodbury formula to iteratively update this choice of $\bfB_k$ as
\begin{equation}\label{eq:Bkupdate}
  \bfB_{k} = \tfrac{k}{k-1}  \bfB_{k-1} \left(\bfI_n - \bfA\t\bfW_k  \left((k-1) \bfI_\ell +  \bfW_k\t\bfA \bfB_{k-1} \bfA\t\bfW_k \right)^{-1} \bfW_k\t\bfA  \bfB_{k-1} \right),
\end{equation}
and for $k = 1$,
\begin{equation}\label{eq:Bkupdate1}
  \bfB_{1} =  \tfrac{1}{\lambda_{1}} \left(\bfI_n - \bfA\t\bfW_1  \left(\lambda_{1}\bfI_\ell +  \bfW_1\t\bfA\bfA\t\bfW_1 \right)^{-1} \bfW_1\t\bfA \right).
\end{equation}
Update formulas~\eqref{eq:Bkupdate} and~\eqref{eq:Bkupdate1} require an inversion of an $\ell \times \ell$ matrix, but $\ell$, which is the number of columns in $\bfW$, is assumed small.

\begin{algorithm}[H] \caption{Stochastic approximation method}
  \begin{algorithmic}[1] \label{alg:SQN}
    \STATE choose initial $\bfx_{0}$, $\lambda_1>0$, suitable sequence $\left\{ \alpha_k \right\}$, and set $k = 1$
    \WHILE{not converged}
      \STATE choose realization $\bfW_k$
      \STATE compute $\bfW_k\t\bfA$ and $\bfW_k\t \bfb$
      \STATE update $\bfB_k$ \label{alg:update}
      \STATE get search direction $\bfs_{k} = - \bfB_k \nabla f_{\bfW_k}(\bfx_{k-1})$
      \STATE update $\bfx_{k} = \bfx_{k-1} + \alpha_k \bfs_{k}$
      \STATE set $k = k+1$
    \ENDWHILE
    \STATE return $ \widehat \bfx = \bfx_{k}$
  \end{algorithmic}
\end{algorithm}

Following most stopping criteria developed for stochastic optimization and stochastic learning methods \cite{Spall2005,Bottou2016,shapiro2014lectures}, we rely on heuristic monitoring of $\bfx_k$ and $f_k = f_{\bfW_k}(\bfx_k)$ to determine early stopping for our methods.  We consider three stopping criteria: a given maximum number of iterations, a certain tolerance on changes in $\bfx_k$, and a tolerance on the improvement in $f_k$, see \cite{Gill1981}.  For the last two, we use
\begin{equation*}
  \norm[\infty]{\bfx_{k-1} - \bfx_k} < \sqrt{{\tt tol}}\left(1+\norm[\infty]{\bfx_k}\right)\\
\quad \mbox{and} \quad
  \abs{\bar f_{k-1} - \bar f_k} < {\tt tol}\left(1+\bar f_{k-1}\right),
\end{equation*}
where ${\tt tol}$ is a given tolerance.
Here, for some number $s\geq 1$, $\bar f_k$ denotes the simple moving average $\bar f_k = mean(f_{k-s},\ldots,f_k)$. In our experiments we choose $s=9$.

Next, we provide two experiments to illustrate our methods.  In Experiment~1, we use a linear regression problem to illustrate various properties or our algorithms and
 validate the theoretical results.  In Experiment~2 we provide some results using a realistic example from machine learning.\\

\subsection*{Experiment 1}
\begin{figure}[bthp]
  \begin{center}
    \includegraphics[width = \textwidth]{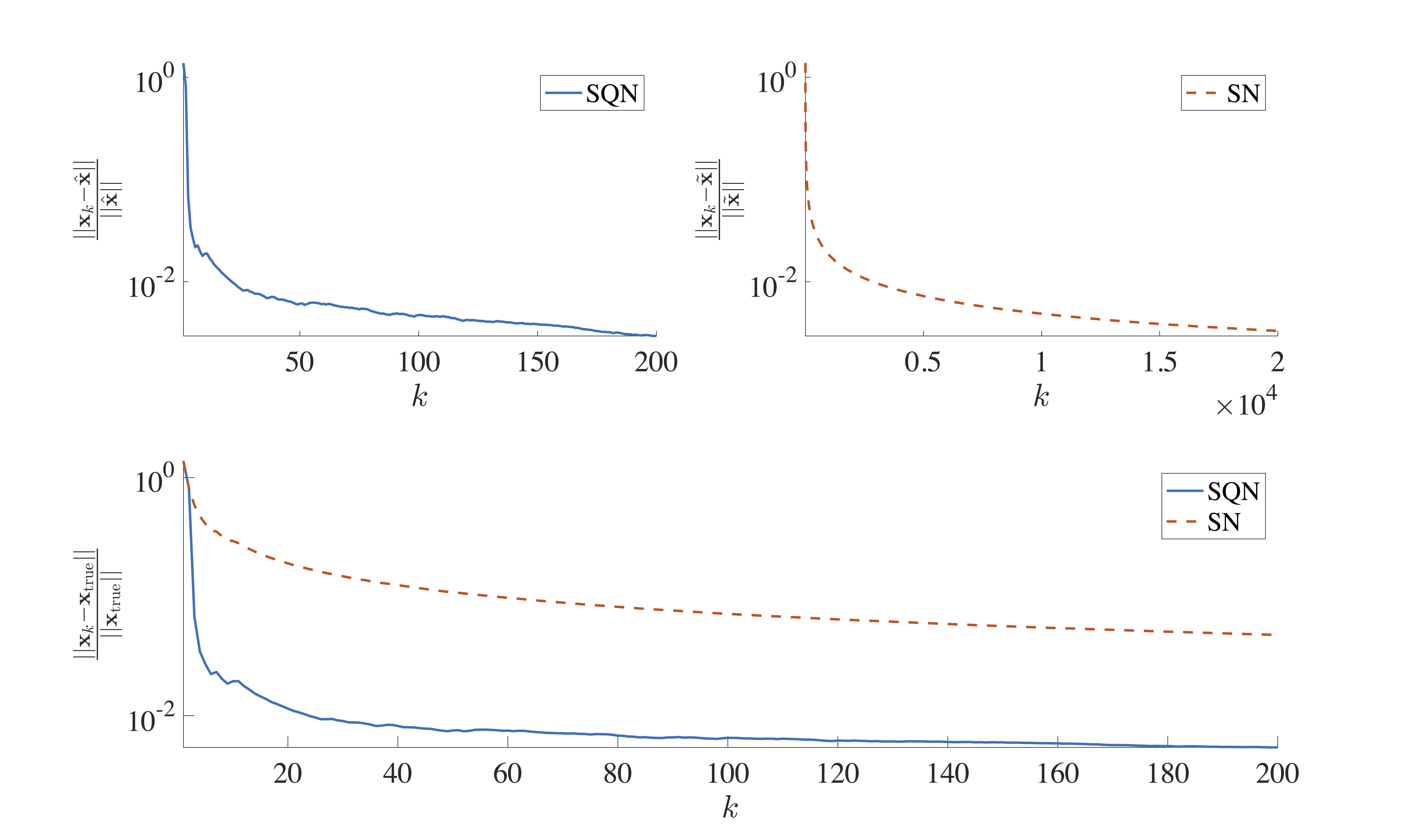}
    \caption{Experiment 1: The top left panel contains relative errors for stochastic quasi-Newton iterates, computed as
$\|\bfx_k - \widehat\bfx\|/\|\widehat \bfx\|$, where $\widehat\bfx$ is the LS solution.  The top right panel contains relative errors for stochastic Newton iterates, computed as
$\|\bfx_k -\widetilde\bfx\|/\|\widetilde \bfx\|$, where $\widetilde\bfx$ is defined in Theorem~\ref{thm:SN}. Notice that we display 20,000 iterations for SN and only 200 iterations for SQN. The bottom panel contains relative errors,
$\|\bfx_k -\bfx_\true\|/\|\bfx_\true\|$, for both SQN and SN.
    } \label{fig:experiment1}
  \end{center}
\end{figure}

We consider a linear regression problem, where $\bfA \in \bbR^{50,000\times 1,000}$ is a random matrix with elements drawn from a standard normal distribution.  We let $\bfx_{\rm true} = \bfone \in \bbR^{1,000}$ and $\bfb = \bfA\bfx_{\rm true} + \bfepsilon$, where the additive noise $\bfepsilon$ is also assumed to be standard normal. We choose $\bfW \in \bbR^{50,000 \times 625}$ to be a block Kaczmarz matrix
with block size $\ell = 625$ and utilize the harmonic sequence
$\alpha_k = {1}/{k}$
as the step size strategy. We select the inverse Hessian approximation $\bfB_k$ using the update formula given in~\eqref{eq:Bkupdate} and~\eqref{eq:Bkupdate1} with $\lambda_1 = 10^{-5}$ 
and start at a random initial guess $\bfx_0$. First, we compare the performance of the stochastic quasi-Newton and Newton methods.
In Figure~\ref{fig:experiment1} we provide plots of relative errors.  In the top left panel, the relative errors are computed as
$\|\bfx_k -\widehat\bfx\|/\|\widehat \bfx\|$, where $\bfx_k$ are stochastic quasi-Newton (SQN) iterates, and in the top right panel, the relative errors are computed as
$\|\bfx_k -\widetilde\bfx\|/\|\widetilde \bfx\|$, where $\bfx_k$ are stochastic Newton  (SN) iterates.  These plots illustrate convergence of SQN and SN iterates to $\widehat\bfx$ and $\widetilde\bfx$ respectively, as proved in Section~\ref{sec:theoretical_results}.
Notice that stochastic Newton (red dashed line) exhibits much slower convergence to $\tilde\bfx$ than the convergence of stochastic quasi-Newton (blue solid line) to $\widehat \bfx$.  In fact, SN requires 20,000 iterations to reach a relative error of $3.3\cdot 10^{-3}$, while the SQN iterates reach a relative error of $3.3\cdot 10^{-3}$ after 175 iterations. Moreover, it takes the stochastic quasi-Newton only 22 iterations to achieve a relative error of $10^{-2}$.


 In the bottom panel of Figure~\ref{fig:experiment1}, we provide reconstruction errors relative to the true solution $\norm[]{\bfx_k -\bfx_\true}/\norm[]{\bfx_\true},$ for both methods, which demonstrates that SQN is faster than SN at providing a better approximation of the true solution. For this experiment the relative error between $\tilde \bfx$ and $\widehat \bfx$ is $\norm[]{\tilde\bfx -\widehat \bfx}/\norm[]{\widehat\bfx} = 5.69\cdot 10^{-3}.$
 We omit the results for the stochastic gradient method due to poor performance.
 It is worth noting that the moderate size of this problem still allows one to
 use a QR solver (e.g., Matlab's ``backslash'') to solve the LS problem, which takes about $6$ seconds whereas SQN requires about $12$ seconds to run $k = 200$ iterations.

\begin{figure}[bthp]
  \begin{center}
    \includegraphics[width = \textwidth]{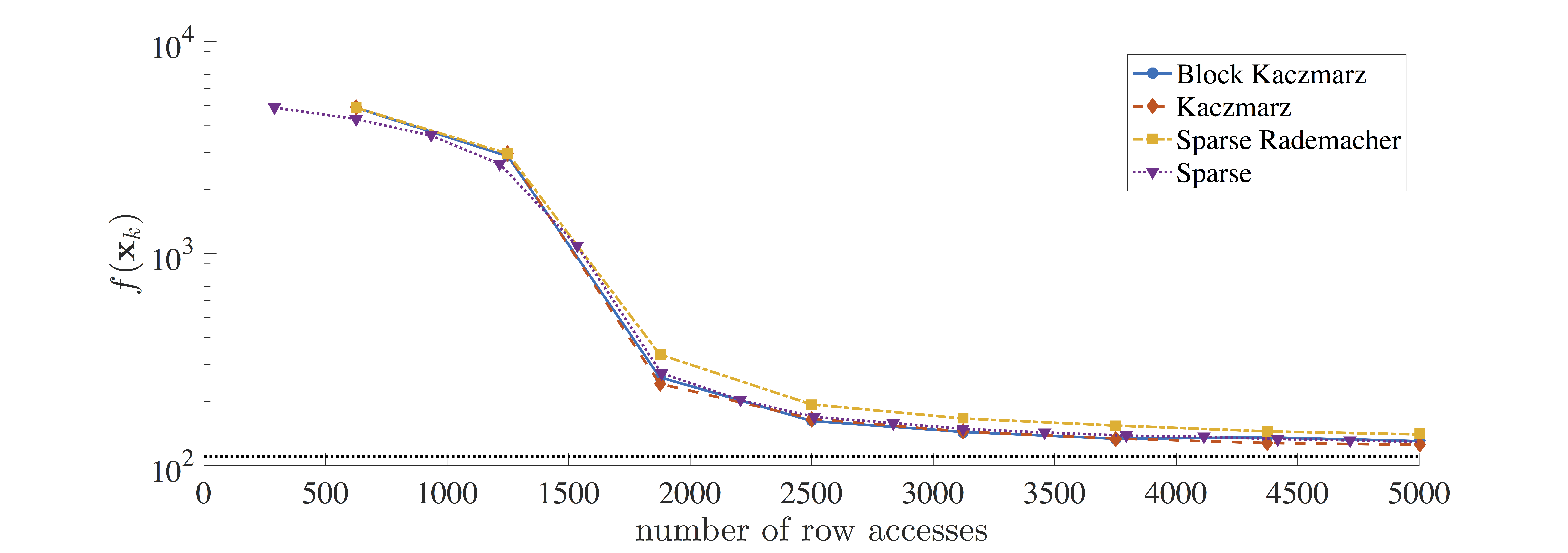}
    \caption{Comparison of function values for different choices of $\bfW$, as a function of the number of row accesses of $\bfA$. Black dotted line corresponds to the LS error $f(\widehat\bfx)$. }   \label{fig:experiment1W}
  \end{center}
\end{figure}

Last, we investigate the performance of four different choices of $\bfW$: (\emph{i}) \emph{Block Kaczmarz}, which is the generalized Kaczmarz method with $\bfQ=\bfI_{m}$ and $\ell_{i}=\ell$ uniformly to create blocks of the same size, (\emph{ii}) \emph{Kaczmarz}, which is similar to the generalized Kaczmarz method with $\bfQ=\bfI_{m},$ but instead of a preset partitioning of the matrix $\bfQ$  into $\bfQ_{i}$ matrices this method samples $\ell$ columns of $\bfQ$ randomly at each iteration to generate $\bfW$, (\emph{iii}) \emph{Sparse Rademacher} with $p=\ell$, and (\emph{iv}) \emph{Sparse Random} with  $\psi = 10^{-5}$.  These choices of $\bfW$ are defined on page~\pageref{pg:kaczmarzW}.
In Figure~\ref{fig:experiment1W}, we provide the function values $f(\bfx_k)$ with respect to the number of row accesses of $\bfA$, since row accesses of $\bfA$ may be the computational bottle neck for large problems. We observe that all choices of $\bfW$ perform similarly, with Kaczmarz and block Kaczmarz having a slight advantage. The performance of the stochastic quasi-Newton method is highly dependent on the underlying problem (e.g., structures in the matrix $\bfA$ and vector $\bfb$) and the realizations of $\bfW$. Empirically, using regression data sets from the UCI Machine Learning Repository \cite{UCIdatasets}, we observe best performances with block Kaczmarz and sparse Rademacher (data not shown).

\subsection*{Experiment 2}
Next, we investigate the use of our stochastic quasi-Newton method for solving large linear LS problems that arise in extreme learning machines (ELMs).  ELM is a machine learning technique that uses random hidden nodes or neurons in a feedforward network to mimic biological learning techniques.  The literature on ELM in the machine learning community is vast, with cited benefits that include higher scalability, less computational complexity, no requirement of tuning, and smaller training errors than generic machine learning techniques. ELM is commonly used for clustering, regression, and classification.  Full details and comparisons are beyond the scope of this paper, and we refer the interested reader to papers such as \cite{huang2011extreme,huang2015trends,huang2012extreme,huang2006extreme,huang2004extreme,huang2015local} and references therein.


At the core of ELM is a very large and potentially dynamically growing linear regression problem.  In this experiment, we investigate the use of stochastic algorithms for efficiently solving these LS problems.
In particular, we consider the problem of handwritten digit classification using the ``MNIST'' database \cite{deng2012mnist}, which contains 60,000 training images and 10,000 testing images of handwritten digits ranging from 0 to 9.  Each image is $28 \times 28$ pixels and converted into a vector $\bfxi \in \bbR^{784}$ (e.g., corresponding to 784 features).

We begin with a brief description of the classification problem for the MNIST dataset.  Suppose we are given a set of $m$ examples in the form of a training set
\begin{equation*}
  S = \left\{(\bfxi_1, c_1), \cdots, (\bfxi_m, c_m) \right\},
\end{equation*}
where $\bfxi_i \in \bbR^{784}$ and $c_i$ takes values from the set of classes $\calC = \left\{0, 1, \cdots, 9 \right\}$.
Consider an ELM with a hidden layer of $n$ nodes, then the goal is to solve a LS problem of the form,
\begin{equation}
  \label{eq:ELMLS}
 \min_\bfX \norm[\fro]{\bfH \bfX - \bfY}^2,
\end{equation}
where the hidden-layer output matrix is defined as
 $$\bfH = \begin{bmatrix} \bfh(\bfxi_1) \\ \vdots \\ \bfh(\bfxi_m)\\ \end{bmatrix}  \in \bbR^{m \times n}$$ with $\bfh(\bfxi) = \begin{bmatrix} h_1(\bfxi) & \cdots & h_n(\bfxi)\end{bmatrix}$ being the output (row) vector of the hidden layer with respect to the input $\bfxi$,
$\bfX=\begin{bmatrix} \bfx_1 & \cdots & \bfx_{10}\end{bmatrix}$ where $\bfx_j\in \bbR^{n}$ contains the desired output weights for class $j$\,,
 and the training data target matrix $\bfY \in \bbR^{m \times 10}$ takes entries
 \begin{equation*}
  y_{ij} = \left\{ \begin{array}{cl} 1, & \mbox{if}\quad c_i = j-1, \\ -1, &  \mbox{else.} \end{array} \right.
\end{equation*}

For this example, $\bfh(\bfxi)$ can be interpreted as a map from the image pixel space to the $n$-dimensional hidden-layer feature space.  Although various activation functions could be used, we employ a standard choice of sigmoid additive hidden nodes with
\[
h_j(\bfxi) = G(\bfd_j, \delta_j, \bfxi) = 1/(1+\exp(-\bfd_j \t \bfxi + \delta_j)),
\]
where all of the hidden-node parameters $(\bfd_j, \delta_j)_{j=1}^n$ are randomly generated based on a uniform distribution \cite{huang2012extreme}.
For our experiments, we set the number of hidden neurons to be $n=300$.

The main computational work of ELM is to solve~\eqref{eq:ELMLS}. Regularized or constrained solutions have been investigated (e.g., \cite{huang2012extreme,luo2014sparse,bai2014sparse}).  However, \emph{our focus will be on solving the unconstrained LS problem efficiently and for very large sets of training data}.  In order to generate larger datasets, we performed multiple random rotations of the original 60,000 training images.  More specifically, each image was rotated by $20(\eta-0.5)$ degrees, where $\eta$ is a random number drawn from a beta distribution with shape parameters equal to $2$.  In our experiments, we consider up to $15$ random rotations per image, resulting in up to 900,000 training images.  Notice that as the number of training images increases, the number of rows of $\bfH$ increases accordingly, while the number of columns remains the same.

We consider three approaches to solve~\eqref{eq:ELMLS} and compare CPU timings.  In the original implementation of ELM \cite{huang2006extreme,huang2004extreme}, the LS solution was computed as $\widehat\bfX = \bfH^\dagger \bfY$ where $\bfH^\dagger$ is the Moore-Penrose  pseudoinverse of $\bfH.$  We denote this approach ``PINV''.  Another approach to solve large (often sparse) LS problems is to use an iterative method such as LSQR \cite{paige1982lsqr,PaSa82b}, but since the LS problem needs to be solved for multiple right-hand sides (here, 10 solves), we use a global least squares method (Gl-LSQR) \cite{toutounian2006global} with a maximum of 50 iterations and a residual tolerance of $10^{-6}$. It was experimentally shown in \cite{toutounian2006global} that Gl-LSQR is more effective and less expensive than LSQR applied to each right hand side.  We use our stochastic quasi-Newton (SQN) method where $\bfW$ corresponds to the sparse Rademacher matrix with $\ell=50$ and $\lambda_1 = 10^{-5}$. We use a maximum number of iterations of 1,000, a stopping tolerance of ${\tt tol} = 10^{-4},$ and an initial guess of $\bfzero$. Since the $\bfB_k$ matrices only depend on $\bfH$ and $\bfW,$ SQN can be applied to multiple right hand sides simultaneously.

\begin{figure}[bthp]
 \begin{center}
   \includegraphics[width=\textwidth]{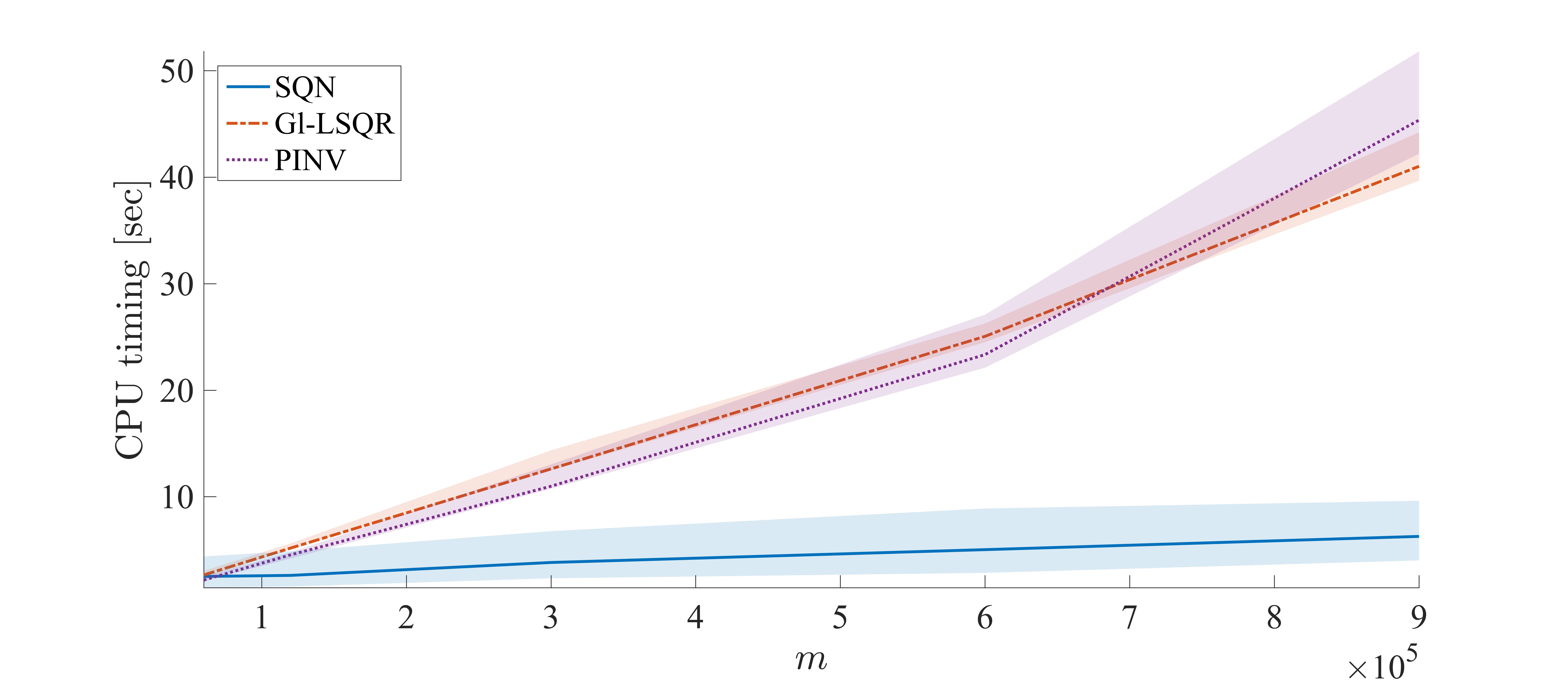}
  \end{center}
  \caption{CPU times (median and 5$^{\rm th}$--95$^{\rm th}$ percentiles) for solving LS problem~\eqref{eq:ELMLS} using stochastic quasi-Newton (SQN), global LSQR (Gl-LSQR), and the Moore-Penrose pseudoinverse for various numbers of training images $m$.}\label{fig:ELMtimings}
\end{figure}

Each LS solver is repeated 20 times in Matlab R2015b on a MacBook Pro with 2.9 GHz Intel Core i7 and 8G memory, and in Figure~\ref{fig:ELMtimings}, we provide the median and 5$^{\rm th}$--95$^{\rm th}$ percentiles of the CPU times vs.\ the number of training images (e.g., number of rows in $\bfH$).  It is evident that for smaller training sets, all three methods perform similarly, but as the number of training images increases, SQN quickly surpasses PINV and Gl-LSQR in terms of faster CPU time.
For various numbers of training data, we provide in Table~\ref{tab:LSerror} the mean and standard deviation of the relative reconstruction error for the SQN estimate, ${\tt rel}=\|\bfX_{\rm SQN}- \widehat \bfX\|_{\rm F}/\|\widehat \bfX\|_{\rm F}$, and of the number of SQN iterations, $k$.  Our results demonstrate that SQN does not necessarily provide the most accurate solutions, however, it can be used to achieve sufficiently good solutions in an efficient manner.

\begin{table}[bthp]
  \caption{For various numbers of training images, we provide the mean and standard deviation for the relative reconstruction errors and the iteration counts for SQN.}
  \label{tab:LSerror}
\begin{center}
  \begin{tabular}{|c|c|c|c|c|c|}
    \hline
    $m$ & $60,\!000$ & $120,\!000$ & $300,\!000$ &$600,\!000$ &$900,\!000$  \\ \hline
    ${\tt rel}$ & 0.2705\,{\small $\pm 0.046$} & 0.2665\,{\small $\pm 0.045$} & 0.2456\,{\small $\pm 0.035$} & 0.2649\,{\small $\pm 0.044$} & 0.2568\,{\small $\pm 0.036$} \\ \hline
    $k$  &  632\,{\small $\pm$219} &  675\,{\small $\pm$227} & 706\,{\small $\pm$194} &  631\,{\small $\pm$196} &
     663\,{\small $\pm$180} \\ \hline
  \end{tabular}
\end{center}
\end{table}

\begin{figure}[bthp]
 \begin{center}
   \includegraphics[width=\textwidth]{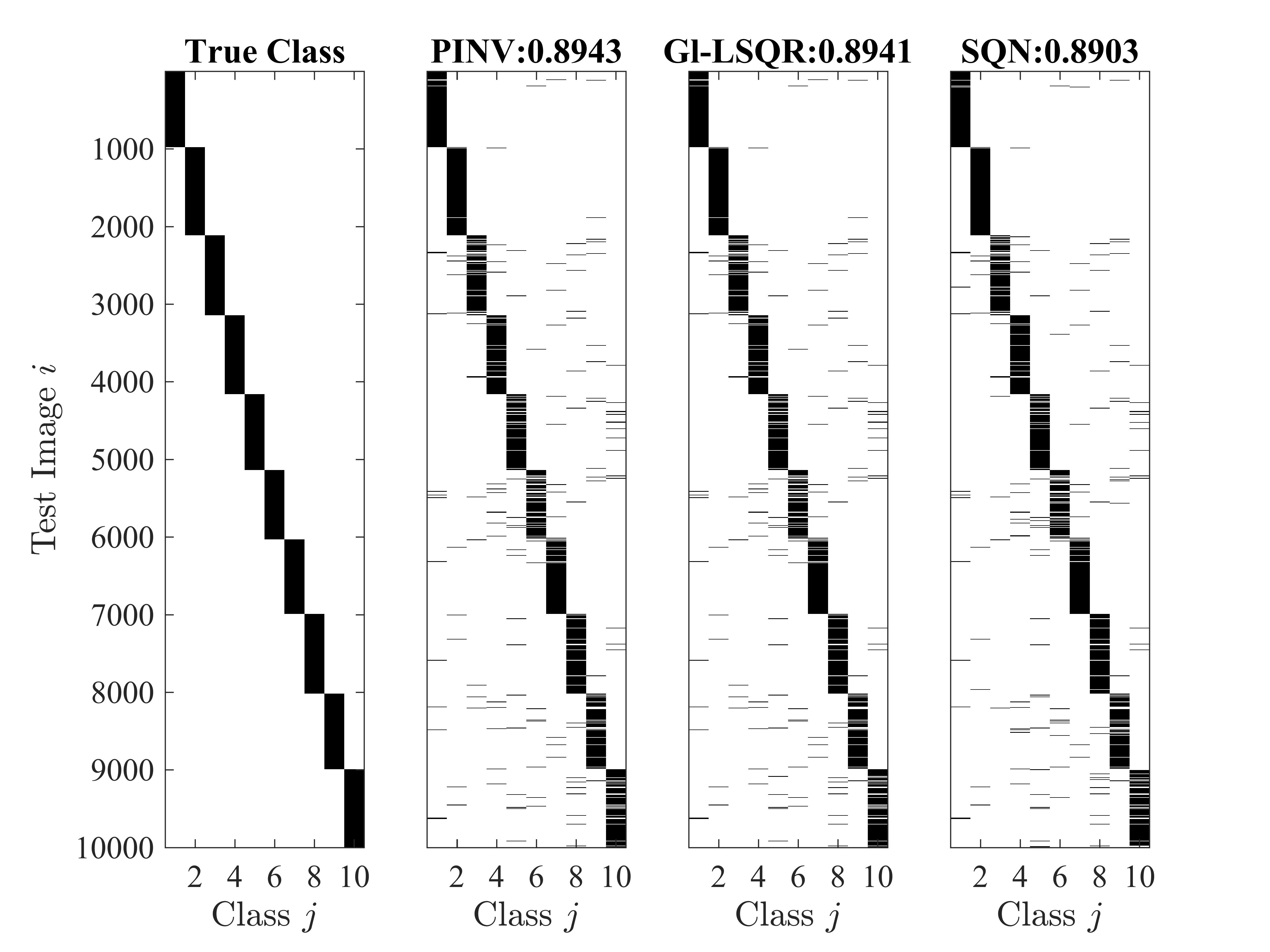}
  \end{center}
  \caption{Classification (with correpsonding accuracy) for the MNIST test images after training on 60,000 images, using different LS solvers.}\label{fig:accuracy}
\end{figure}

Next we test the performance of these estimates for classification of the MNIST testing dataset.  That is, once computed, the output weights $\bfX$ can be used to classify images in the following way. For each test image, the predicted class is given by
\begin{equation*}
 \mbox{Class of } \bfxi = \arg\max_j \bfh(\bfxi) \bfx_j \,.
\end{equation*}
In Figure~\ref{fig:accuracy} we provide a visualization of the computed classifications for the 10,000 testing images, where accuracy values in the titles are calculated as $1-{r}/{10000}$ where $r$ is the number of misclassified images.  An accuracy value that is close to $1$ corresponds to a good performance of the classifier.  These results correspond to training on 60,000 images, and the testing set was sorted by class for easier visualization.
Notice that in Figure~\ref{fig:accuracy} the misclassified images are almost identical for all three methods, and the classification accuracy for SQN is only slightly smaller than that of PINV and Gl-LSQR.  Thus, we have shown that our SQN method can achieve comparable classification performance as PINV and GL-LSQR with much faster learning speed.

It is worth noting that the matrices considered here, though large, can still be loaded into memory. For problems where this is not the case (e.g., data too large or being dynamically generated \cite{zhao2012online}), PINV and Gl-LSQR would not be feasible, while SQN could still be used.

\section{Conclusions} 
\label{sec:conclusions}
In this paper we introduced a stochastic Newton and a stochastic quasi-Newton method for solving very large linear least-squares problems.
New theoretical results show that under mild regularity conditions on the distributions, the stochastic quasi-Newton iterates converge to the unique LS solution
$\widehat\bfx$, while the stochastic Newton iterates converge to a different estimator of $\bfx_\mathrm{true}$ which is still unbiased but is likely to
have larger variance.
We provide implementation details for choices of random variables $\bfW_k$, step length parameters $\alpha_k$, and stochastic inverse Hessian approximations $\bfB_k$. Our numerical experiments validate our theory and also show the potential benefits of these methods for machine learning applications where the data sets of interests are typically extremely large.  Furthermore, by allowing sparse distributions, our approaches can handle scenarios where the entire matrix $\bfA$ is not available or full matrix-vector multiplications with $\bfA$ are not feasible.

This work opens the doors to a plethora of new methods, advancements, and implementations for stochastic Newton and quasi-Newton methods.  Current work includes the development of limited-memory quasi-Newton methods to reduce the storage and updating costs for the sequence of $\bfB_{k}$ matrices, adaptive choices of $\bfW$ that can exploit matrix properties, efficient parallel implementations, and extensions to nonlinear problems.

\bigskip
\acknowledgement{The authors are grateful to Qi Long for initial discussions and Arvind Saibaba for helpful comments that have improved the manuscript. We are also grateful to Eldad Haber for pointing us to the ELM problem and to the Institute for Mathematics and its Applications for nurturing this research collaboration.}


\section{Appendix}\label{sec:appendix}

\subsection{Proof of Theorem \ref{thm:SQN}}\label{sec:appendixSQN}
\emph{(i)}  Let $\bfv \in \mathbb{R}^{m}$ and assume that $\bfW$ has the same distribution as $\bfW_k$. Then,
		\[
			\|\bfW\bfW^\top \bfv\|^2 \leq \norm[]{\bfW\bfW^\top}^2\|\bfv\|^2,
		\]
		and therefore  $\Ex(\,\bfv^\top \bfW\bfW^\top\bfW\bfW^\top\bfv\,) < \infty$. Furthermore,
		\[
			\bfv^\top\bfC\bfv = \Ex(\,\bfv^\top \bfW\bfW^\top\bfW\bfW^\top\bfv\,)
			 = \Ex(\,\|\bfW\bfW^\top \bfv\|^2\,) \geq 0,
		\]
		and therefore $\bfv^\top\bfC\bfv=0$ iff $\|\bfW\bfW^\top \bfv\|^2=0$ a.s., which happens iff $\bfW\bfW^\top \bfv=\zero$ a.s., in which case
		\[
			\Ex(\,\bfW\bfW^\top\,)\, \bfv = \beta\,\bfI_{n}\bfv = \zero.
		\]
		This implies $\bfv=\zero$, and thus  $\bfC$ is defined and SPD.

\emph{(ii)} Note that
		\begin{equation}\label{eq:Bskr}
			\bfs_k = -\bfB_k\bfA^\top\bfW_k\bfW_k^\top\bfA(\bfx_{k-1}-\widehat{\bfx})   -\bfB_k\bfA^\top\bfW_k\bfW_k^\top\bfr,
		\end{equation}
		where
		\(
			\bfr = \widehat{\bfb} - \bfb.
		\)
		Therefore,
		\[
			\bfb_{k+1} - \widehat{\bfb} = \bfb_k - \widehat{\bfb} - \alpha_k\,\bfA\bfB_k\bfA^\top\bfW_k\bfW_k^\top (\bfb_k-\widehat{\bfb}) -\alpha_k\,\bfA\bfB_k\bfA^\top\bfW_k\bfW_k^\top\bfr,
		\]
		and there are positive constants $a, b_k$ and $c_k$ such that
		\[
			e_{k+1} \leq  a\,e_k + b_k \,B^2 \norm[]{\bfW_k\bfW_k^\top}^2\,e_{k} + c_k \,B^2 \norm[]{\bfW_k\bfW_k^\top}^2
		\]
		for all $k$.  Since $\bfW_k$ and $e_k$ are independent and $\Ex\,e_0 < \infty$, it follows that $\Ex\,e_k<\infty$ for all $k$.
		This together with \eqref{eq:Bskr} implies that $\Ex(\,\norm[]{\bfs_k}^2 \,)<\infty$ for all $k$.

\emph{(iii)} By definition,
		\[
			e_{k+1} - e_k = 2\alpha_k\,(\bfb_k-\widehat{\bfb})^\top\bfA\bfs_k + \alpha_k^2\, \bfs_k^\top\bfA^\top\bfA\bfs_k,
		\]
		and therefore
		\[
			\Ex(\,e_{k+1} - e_k \mid \Fc_k\,) = 2\alpha_k\,(\bfb_k-\widehat{\bfb})^\top\bfA\,\Ex(\,\bfs_k\mid \Fc_k\,)  + \alpha_k^2\, \Ex(\,\bfs_k^\top\bfA^\top\bfA\bfs_k\mid \Fc_k\,).
		\]
		Since this residual vector is orthogonal to the column space of $\bfA$, \eqref{eq:Bskr}  leads to
		\begin{equation} \label{eq:EskB}
				\Ex (\,\bfs_{k}\mid\Fc_k \,)  = -\bfB_k\bfA^\top \Ex(\bfW_k\bfW_k)(\bfA\bfx_{k-1} -\bfA\widehat{\bfx}) = -\beta\,\bfB_k\bfA^\top(\bfb_k-\widehat{\bfb}).
		\end{equation}
		Let $\lambda_{\bfA}^2:=\lambda_\mathrm{max}(\bfA\bfA^\top)$, then
		\begin{align*}
			\bfs_k^\top\bfA^\top\bfA\bfs_k &= (\bfb_k-\bfb)^\top\bfW_k\bfW_k^\top\bfA\bfB_k\bfB_k\bfA^\top\bfW_k\bfW_k^\top(\bfb_k-\bfb)\\
			& \leq  \lambda_{\bfA}^2 B^2\,(\bfb_k-\bfb)^\top\bfW_k\bfW_k^\top\bfW_k\bfW_k^\top(\bfb_k-\bfb),
		\end{align*}
		and therefore
		\begin{equation} \label{eq:EsAAs}
			\Ex(\,\bfs_k^\top\bfA^\top\bfA\bfs_k \mid \Fc_k\,) \leq \lambda_{\bfA}^2 B^2\,\norm[\bfC]{\bfb_k-\bfb}^2\leq 2 \lambda_{\bfA}^2 B^2\,e_k + 2 \lambda_{\bfA}^2 B^2\,\norm[\bfC]{\bfr}^2.
		\end{equation}
		Equations~\eqref{eq:EskB} and~\eqref{eq:EsAAs} lead to
		\[
			\Ex(\,e_{k+1}-e_k\mid \Fc_k\,) \leq  -2\beta\,\alpha_k\,(\bfb_k-\widehat{\bfb})^\top\bfA\bfB_k\bfA^\top(\bfb_k-\widehat{\bfb}) +  2 \alpha_k^2\lambda_{\bfA}^2 B^2\,e_k + 2\alpha_k^2\lambda_{\bfA}^2 B^2 \norm[\bfC]{\bfr}^2,
		\]
		subtracting $2\alpha_{k}^{2}\lambda_{\bfA}^{2}B$ from both sides yields
		\[
			\Ex(\,e_{k+1}-(1+2 B^2 \lambda_{\bfA}^2\alpha_k^2)e_k\mid \Fc_k\,) \leq  -2\beta\,\alpha_k\,(\bfb_k-\widehat{\bfb})^\top\bfA\bfB_k\bfA^\top(\bfb_k-\widehat{\bfb})  + 2\lambda_{\bfA}^2 \alpha_k^2\,B^2 \norm[\bfC]{\bfr}^2.
		\]
		Define
		\(
			\gamma_k = \prod_{i=1}^{k-1}(\,1+2B^2\lambda_{\bfA}^2\alpha_k^2\,)^{-1}\leq 1.
		\)
		Since
		\[
			0\leq -\log(\gamma_k) \leq 2 B^2 \lambda_{\bfA}^2\sum_{i=1}^{\infty}\alpha_i^2 <\infty
		\]
		for any $k$, it follows that $\{\gamma_{k}\}$ is a decreasing and convergence sequence to  some $\gamma>0$. Define $\tilde{e}_k = \gamma_ke_k$ and $Z_k=\bbE(\,\tilde{e}_{k+1} - \tilde{e}_{k} | \calF_{k} \,)$,
 then
		\begin{eqnarray}
			\Ex(\tilde{e}_{k+1} -\tilde{e}_{k}\mid\Fc_k)
			& \leq &- 2\beta\alpha_k\gamma_{k+1} (\bfb_k-\widehat{\bfb})^\top\bfA\bfB_k\bfA^\top(\bfb_k-\widehat{\bfb}) + 2 B^2\gamma_{k+1} \alpha_k^2 \lambda_{\bfA}^2 \norm[\bfC]{\bfr}^2 \nonumber\\
			& \leq & 2 B^2 \alpha_k^2\, \lambda_{\bfA}^2 \norm[\bfC]{\bfr}^2, \label{eq:Bconddiff1}
		\end{eqnarray}
%
		and also
		\[
			\sum_{k=1}^{\infty}\Ex[\,\I_{Z_k\geq 0}\,(\tilde{e}_{k+1}-\tilde{e}_k)\,]\,\, \leq \,\,2 B^2 \lambda_{\bfA}^2\norm[\bfC]{\bfr}^2\sum_{k=1}^{\infty}\alpha_k^2<\infty.
		\]
		Theorem~\ref{thm:quasi} now implies that $\{\tilde{e}_k\}$ converges a.s.\ and since $\gamma_k$ converges to a nonzero value, it follows that $\{e_k\}$ also converges a.s. The final step is to show that $e_k\arras 0$. It follows from \eqref{eq:Bconddiff1} that
		\begin{align*}
			2\beta \sum_{k=1}^n \alpha_k\gamma_{k+1}\,\Ex&(\,\lambda_{\min}(\bfB_k)\|\,\bfA^\top(\bfb_k-\widehat{\bfb})\,\|^2\,) \\
			&\leq  2\beta\,\sum_{k=1}^n  \alpha_k\gamma_{k+1}(\bfb_k-\widehat{\bfb})^\top\bfA\bfB_k\bfA^\top(\bfb_k-\widehat{\bfb})\\
			& \leq 2 B^2 \lambda_{\bfA}^2\|\bfr\|_{\bfC}^2\sum_{k=1}^n\alpha_k^2  + \Ex\,\tilde{e}_1 <\infty
		\end{align*}
		for any $n$. Therefore,
		\[
		\gamma\sum_{k=1}^\infty \alpha_k\, \lambda_{\min} (\bfB_k)\,\|\bfA^\top(\bfb_k-\widehat{\bfb})\|^2
			 \leq \sum_{k=1}^\infty \alpha_k\gamma_{k+1}\,\lambda_{\min}(\bfB_k)\,\|\bfA^\top(\bfb_k-\widehat{\bfb})\|^2 < \infty \quad \mathrm{a.s.}
                 \]
%
		Since $\sum_{k=1}^{\infty} \alpha_k \lambda_{\min}(\bfB_k)= \infty$ a.s., it follows that
		for almost all $w\in \Omega$  there is a subsequence $n_k(w)$ such that
		\(
			\norm[(\bfA^\top\bfA)^2]{\bfx_{n_k(w)}(w)-\widehat{\bfx}}^2 \,\, \to\,\, 0,
		\)
		which also implies
		\[
			\|\,\bfx_{n_k(w)}(w)-\widehat{\bfx}\,\|^2= e_{n_k(w)}(w) \,\, \to\,\, 0.
		\]
		and therefore, since $e_k$ converges a.s., we also have $e_k\arras 0$. Hence, $\bfx_k\arras \widehat{\bfx}$.

\subsection{Proof of Theorem \ref{thm:SN}} \label{sec:appendixSN}

	The proof is a slight modification of that of Theorem~\ref{thm:SQN}. Define $\Fc_k = \sigma(\,\bfW_i;\,i<k\,)$.
		\emph{(i)} Let $\bfv \in \mathbb{R}^{m}.$ Since
		\(
			\bfv\t\bfC\bfv = \bbE(\,\|\bfH_k\bfv\|^{2}\,)  \geq 0,
		\)
		it follows that $\bfC$ is semi-positive definite, and therefore so is $\bfA\t\bfC\bfA$. Since $(\bfW_k\t\bfA)^\dagger\bfW_k\t\bfA$ is symmetric, $\bfP\bfA$ is symmetric.
 Let $\bfv \in\bbR^{n}.$ Using properties of the pseudoinverse  gives
		\[
		\bfv\t\bfP\bfA\bfv=	\bfv\t\bbE(\,\bfH_k\bfA\,)\bfv = \bbE (\,\|\,\bfH_k\bfA\bfv\,\|^{2}  \,)  \geq 0,
 		\]
		where equality holds iff $\bfH_k\bfA\bfv=\zero$ a.s., and since $\Ex(\bfW_k\bfW_k\t)=\beta\bfI_m$ and $\bfA$ is full column-rank, it follows that $\bfv=\zero$. Hence
$\bfP\bfA$ is SPD.

		\emph{(ii)} Note that
		\begin{equation}\label{eq:Bskrr}
			\bfs_k = -\bfH_k\bfA(\bfx_{k-1}-\widetilde{\bfx})   -\bfH_k\,\bfr,
		\end{equation}
where $\bfr = \bfA\widetilde{\bfx}-\bfb$. Therefore
		\[
			\bfx_{k}- \widetilde{\bfx} = \bfx_{k-1} - \widetilde{\bfx}- \alpha_{k}\,\bfH_k\bfA(\bfx_{k-1} - \widetilde{\bfx}) - \alpha_{k}\,\bfH_{k}\,\bfr,
		\]
and since $\bfH_k\bfA$ is an orthogonal projection matrix, it follows that
\[
	e_{k+1} \leq  4(1+\alpha_k^2)\,e_k +  4\alpha_k^2 \,\|\bfH_k\bfr\|^2.
\]
This then leads to
\[
\Ex\,e_{k+1} \leq  4(1+\alpha_k^2)\,\Ex\,e_k +  4\alpha_k^2 \,\|\bfC\bfr\|^2,
\]
which implies that $\Ex\,\e_k<\infty$ and $\Ex(\,\|\bfs_k\|^2\,)<\infty$ for all $k$.

 		\emph{(iii)} The idea is again to use the quasimartingale convergence theorem. Since $\bfW_k$ and $\Fc_k$ are independent,
\begin{eqnarray}
\bbE(\, e_{k+1}-e_{k} \,\mid\, \calF_{k} \,) &= & 2\alpha_{k}\,(\bfx_{k-1}-\widetilde{\bfx})\t\bbE(\,\bfs_{k} \, \mid \, \calF_{k}\,) + \alpha_{k}^{2}\, \bbE(\,\|\bfs_{k}\|^{2} \, \mid \, \calF_{k}\,)\nonumber \\
&=& - 2\alpha_{k}\norm[\bfP\bfA]{\bfx_{k-1} - \widetilde{\bfx}}^{2} + \alpha_{k}^{2}\, \bbE(\,\|\bfs_{k}\|^{2} \, \mid \, \calF_{k}\,).\label{eq:firste1}
\end{eqnarray}
To put a bound on the second term of~\eqref{eq:firste1} we use again the fact that $\bfH_{k}\bfA$ is a projection matrix to obtain:
\begin{equation} \label{eq:eigen}
 \bbE(\,\|\bfs_{k}\|^2 \mid  \calF_{k} \,)  \leq  2(\bfA\widetilde{\bfx}-\bfb)\t\bfC(\bfA\widetilde{\bfx}-\bfb) + 2\|\bfx_{k-1}-\widetilde{\bfx}\|^{2}
				 \leq  c_{1} + 2 \|\bfx_{k-1}-\widetilde{\bfx}\|^{2},
 \end{equation}
where $c_{1} =\lambda_{\text{max}}\left(\bfC\right)\norm[]{\bfA\widetilde{\bfx}-\bfb}^{2}$. Therefore, equation  \eqref{eq:firste1} can be bounded as
 	 	\begin{equation*} \label{eq:firste2}
	 		\bbE(\, e_{k+1}-e_{k} \mid \calF_{k}\,) \leq -2\alpha_{k}\norm[\bfP\bfA]{\bfx_{k-1}-\widetilde{\bfx}}^{2}+ c_1 \alpha_{k}^{2}  + 2\alpha_{k}^{2}\,e_{k},
	 	\end{equation*}
which yields
\begin{equation} \label{eq:bound1}
 \bbE(\, e_{k+1}- e_{k}(1+ \alpha_{k}^{2}) \mid  \calF_{k} \,) \leq -2\alpha_{k}\,\norm[\bfP\bfA]{\bfx_{k-1}-\widetilde{\bfx}}^{2}+ c_1 \alpha_{k}^{2}  \leq c_1 \alpha_{k}^{2}.
\end{equation}
Let
\(
\nu_{k} = \prod_{i=1}^{k-1}(\,1+\alpha_{i}^{2}\,)^{-1}<1.
\)
The sequence $\{\nu_{k}\}$ converges to some $\nu>0$ because $\sum_{k=1}^{\infty} \alpha_{k}^{2} < \infty$. Define $\widetilde{e}_{k}=\nu_{k}e_{k},$ and multiply both sides of \eqref{eq:bound1} by $\nu_{k+1}$. We obtain
\begin{equation} \label{eq:rearg}
	\bbE(\,\widetilde{e}_{k+1}-\widetilde{e}_{k} \mid \calF_{k}\,) \leq  -2\alpha_{k}\,\nu_{k+1}\norm[\bfP\bfA]{\bfx_{k-1}-\widetilde{\bfx}}^{2} +  c_1\alpha_{k}^{2}\,\nu_{k+1}  \leq  c_1\alpha_{k}^{2}\,\nu_{k+1} .
\end{equation}
Define $Z_k = \bbE(\,\widetilde{e}_{k+1}-\widetilde{e}_{k} \mid \calF_{k}\,)$. Then,
\[
\bbE [\,\I_{Z_{k}\geq 0}\,(\widetilde{e}_{k+1}-\widetilde{e}_{k})\,]  = \bbE(\,\I_{Z_k\geq 0}\,\bbE [\,\widetilde{e}_{k+1}-\widetilde{e}_{k} \mid \calF_{k}\,]\,)  \leq c_1\alpha_{k}^{2}\,\nu_{k+1}.
\]
Since $\sum_{k=1}^{\infty} \alpha_{k}^{2}\nu_{k+1} < \infty$ the series  $\sum_{k=1}^{\infty} \bbE (\, \I_{X_k\geq 0}\,(\tilde{e}_{k+1} - \tilde{e}_{k})\,)$ converges and therefore $\{\tilde{e}_k\}$ converges a.s. by Theorem~\ref{thm:quasi}. But since $\nu_k$ converges to a nonzero value, it also follows that $\{e_k\}$ converges a.s. The final step is to show that $\{e_k\}$ in fact converges to zero. Rearranging the terms and taking the expected value of both sides of \eqref{eq:rearg} yields
\[
\sum_{k=1}^{\infty}\alpha_{k}\nu_{k+1}\bbE(\,\norm[\bfP\bfA]{\bfx_{k-1}-\widetilde{\bfx}}^2\,) < \infty,
\]
Since $\sum_{k=1}^{\infty} \alpha_{k}=\infty$ and $\nu_{k}\to\nu>0$, it follows that $\bbE(\,\norm[\bfP\bfA]{\bfx_{n_k}-\widetilde{\bfx}}^2\,)\to 0$ for some subsequence $(n_k)$, and therefore we also have $\bbE(\,\|\bfx_{n_k}-\widetilde{\bfx}\|^2\,)=\Ex\,e_{n_k}\to 0$.  By Fatou's lemma:
\[
0 \leq \Ex\lim\inf e_{n_k} =\Ex\lim e_{n_k} \leq \lim\inf \Ex\,e_{n_k} = 0.
\]
It follows that $\lim e_{n_k}=0$ a.s.\ and since $\{e_k\}$ converges a.s., this implies $e_k\arras 0$ and therefore $\bfx_k\arras \widetilde \bfx$.

\bibliographystyle{abbrv}
\bibliography{stochref}
\end{document}